%% file: CMBE24Abstract.tex
\documentclass{cmbe24}
\usepackage{amsmath,amsfonts,amssymb}
\usepackage{graphicx,wrapfig}
\usepackage{hyperref}
\usepackage{import}
\usepackage{xcolor}
\usepackage{subfigure}
\usepackage{wrapfig}
\usepackage{siunitx}
\usepackage{pgfplots}
\pgfplotsset{compat=newest}
\DeclareSIUnit{\mmHg}{mmHg}
\usepackage{lipsum}

\def\vct{\mathbf}

\title{A coupled fluid-dynamics-heat transfer model for 3D simulations of the aqueous humor flow in the human eye}

\author[1]{Thomas Saigre}
\author[1]{Christophe Prud'homme}
\author[2]{Marcela Szopos}
\author[1]{Vincent Chabannes}
\affil[1]{Institut de Recherche Mathématique Avancée, UMR 7501 Université de Strasbourg et CNRS \iffalse,\texttt{thomas.saigre@math.unistra.fr, christophe.prudhomme@cemosis.fr, vincent.chabannes@cemosis.fr}\fi}
\affil[2]{Université Paris Cité, CNRS, MAP5, F-75006 Paris, France \iffalse, \texttt{marcela.szopos@u-paris.fr}\fi}

\summary{Understanding human eye behavior involves intricate interactions between physical phenomena such as heat transfer and fluid dynamics. Accurate computational models are vital for comprehending ocular diseases and therapeutic interventions.
This work focuses on modeling and simulating aqueous humor flow in the anterior and posterior chambers of the eye, coupled with overall heat transfer.
Aqueous humor dynamics regulates intraocular pressure, crucial for understanding conditions like glaucoma.
Convective effects from temperature disparities also influence this flow.
Extending prior research, this work develops a comprehensive three-dimensional computational model to simulate coupled fluid-dynamic-heat transfer model, thus contributing to the understanding of ocular physiology.
}

\keywords{mathematical and computational ophthalmology, finite element method, thermo-fluid dynamics}

\begin{document}

% section without full stop
\section{Introduction}

Understanding the behavior of the human eye is a challenging task, as it involves the study of the interaction between different physical phenomena, such as heat transfer, fluid dynamics, and tissue deformation.
It is therefore crucial to develop accurate and efficient computational models to simulate complex multi-physics physiology of the ocular system, in order to gain a better understanding of the mechanisms of ocular diseases and the associated therapeutical interventions.
In the present work, we focus on the simulation of the flow of the aqueous humor (AH) in the anterior and posterior chambers of the human eye, and its coupling with the heat transfer inside the eyeball.%, \emph{using pressure boundary conditions}\textbf{à compléter ou enelver}.

The AH is a transparent fluid produced by the ciliary body (see Fig.~\ref{fig:geo-eye}), that flows from the posterior chamber (PC) to the anterior chamber (AC), where it is drained through two pathways, the trabecular meshwork, and the uveoscleral pathway.
AH flow plays a fundamental role in the maintenance of the intraocular pressure (IOP) level.
IOP is a key parameter since when abnormally elevated, it is a major risk factor for degenerative ocular diseases such as glaucoma~\cite{glaucoma}.
In addition to the hydraulic pressure difference created by production and drainage, AH dynamic is influenced by posture and thermal factors.
Specifically, there are convective effects produced by the temperature difference between the external environment at the corneal surface and the internal surface, which is at the body temperature.

Numerical investigations of this complex dynamics have been proposed by several authors, such as \cite{Wang,Murgoitio} where the model describes flow coupled with heat transfer solely in the AC and PC;
\cite{PoF} where the impact of pressure on the AH flow and drainage is studied; or \cite{Ooi,Abd} where a coupled model for the thermo-fluid-dynamics of the AH flow in the AC with boundary conditions on the velocity is investigated.
The present work develops a three-dimensional computational model allowing to simulate the heat transfer inside the whole human eyeball, taking into account the flow of AH in the AC and PC.

This contribution represents an extension of the work presented in \cite{Heat}, where heat transfer inside the human eyeball was investigated with a special focus on the parametric effects.

\section{Methodology}

\begin{figure}
    \centering
    \subfigure[Geometrical model of the human eye.]{
        \label{fig:geo-eye}
        \def\svgwidth{0.5\columnwidth}
        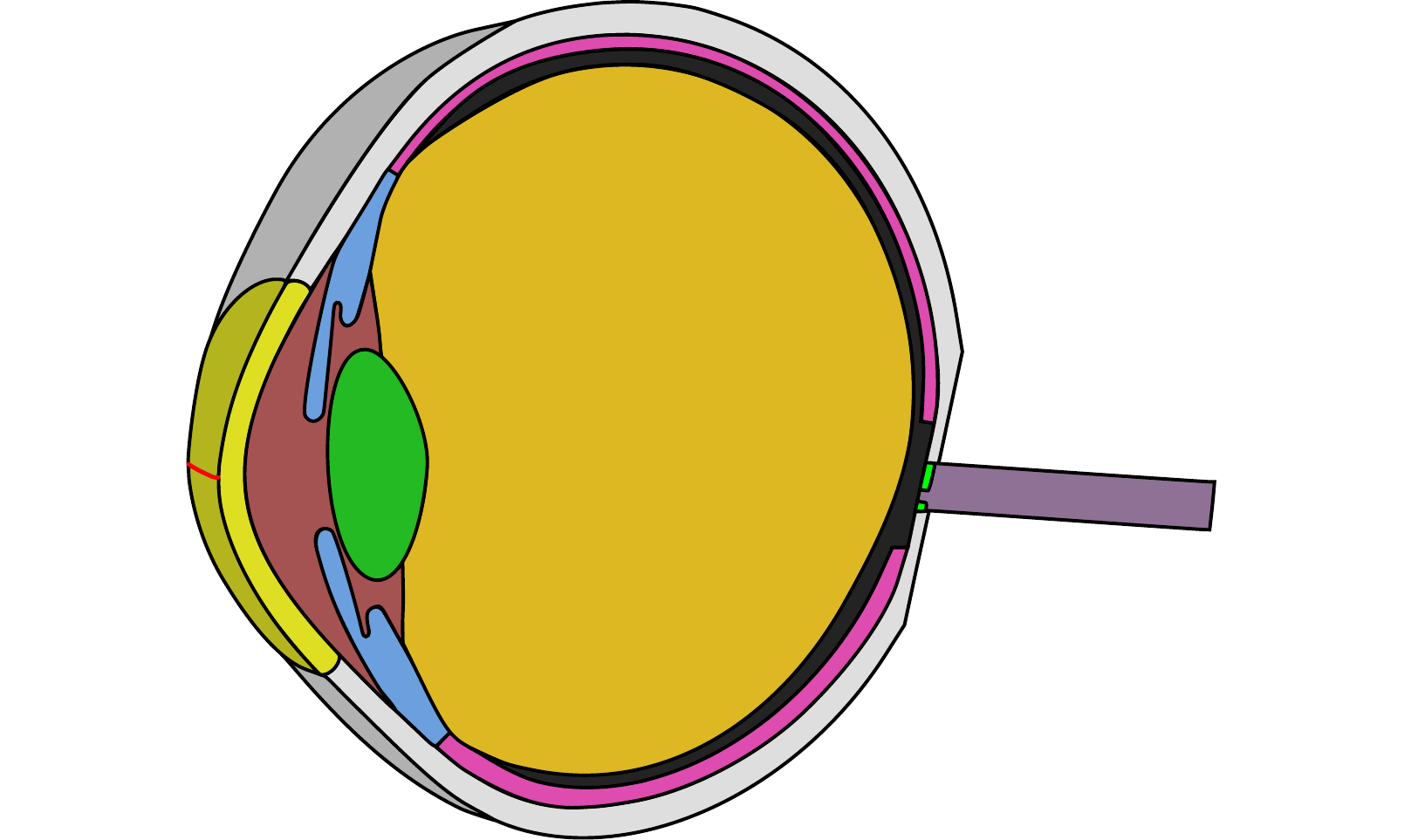
    }
    \subfigure[Vertical cut of the anterior and posterior chambers of the eye.]{
        \label{fig:ah}
        \def\svgwidth{0.3\columnwidth}
        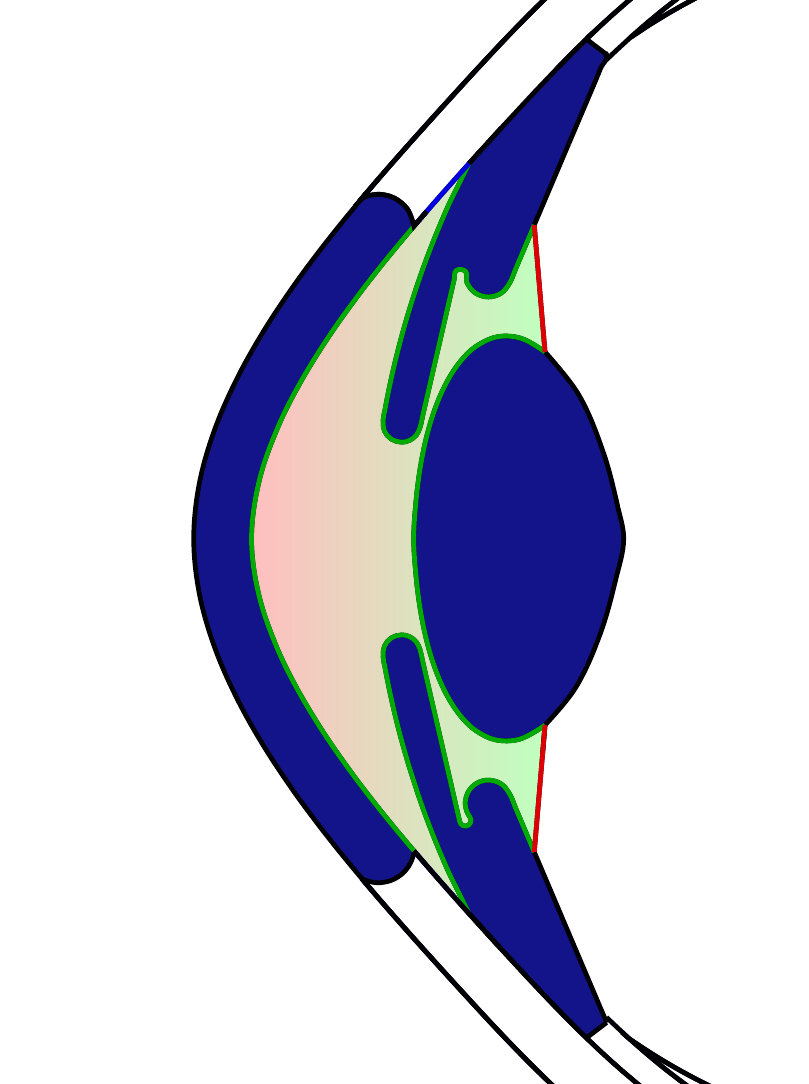
    }
    \caption{Description of the human eyeball, and the anterior and posterior chambers of the eye.}
\end{figure}

We denote by $\Omega$ the domain of the human eye presented in Fig.~\ref{fig:geo-eye}.
This domain is divided in ten subdomains, representing the different parts of the eye, such as the cornea, the lens, the vitreous body, the retina, \emph{etc.}, with various physical properties.
Precisely, we denote by $\Omega_\text{AH}$ the anterior and posterior chambers of the eye (brown part in Fig.~\ref{fig:geo-eye}),
which is filled with the AH, as shown in Fig.~\ref{fig:ah}.

\begin{wrapfigure}{L}{0pt}
    \centering
    \begin{tikzpicture}
    \begin{axis}[
        colorbar,
        colormap/jet, % Choose the colormap you prefer
        % axis equal image,
        enlargelimits=false,
        point meta max=310.14813232421875,
        point meta min=308.0091857910156,
        colorbar horizontal,
        axis line style = {draw=none},
        tick style = {draw=none},
        xtick = \empty, ytick = \empty,
        colorbar style={
            xlabel = {$T$ [\si{\kelvin}]},
            height=0.05*\pgfkeysvalueof{/pgfplots/parent axis height},
            width=0.9*\pgfkeysvalueof{/pgfplots/parent axis width},
            at={(0.5,-0.02)}, % Adjust the position to center vertically
            anchor=center, % Adjust the anchor point
        },
        width=0.4\textwidth
    ]
        \addplot graphics [includegraphics cmd=\pgfimage, xmin=0, xmax=1, ymin=0, ymax=1] {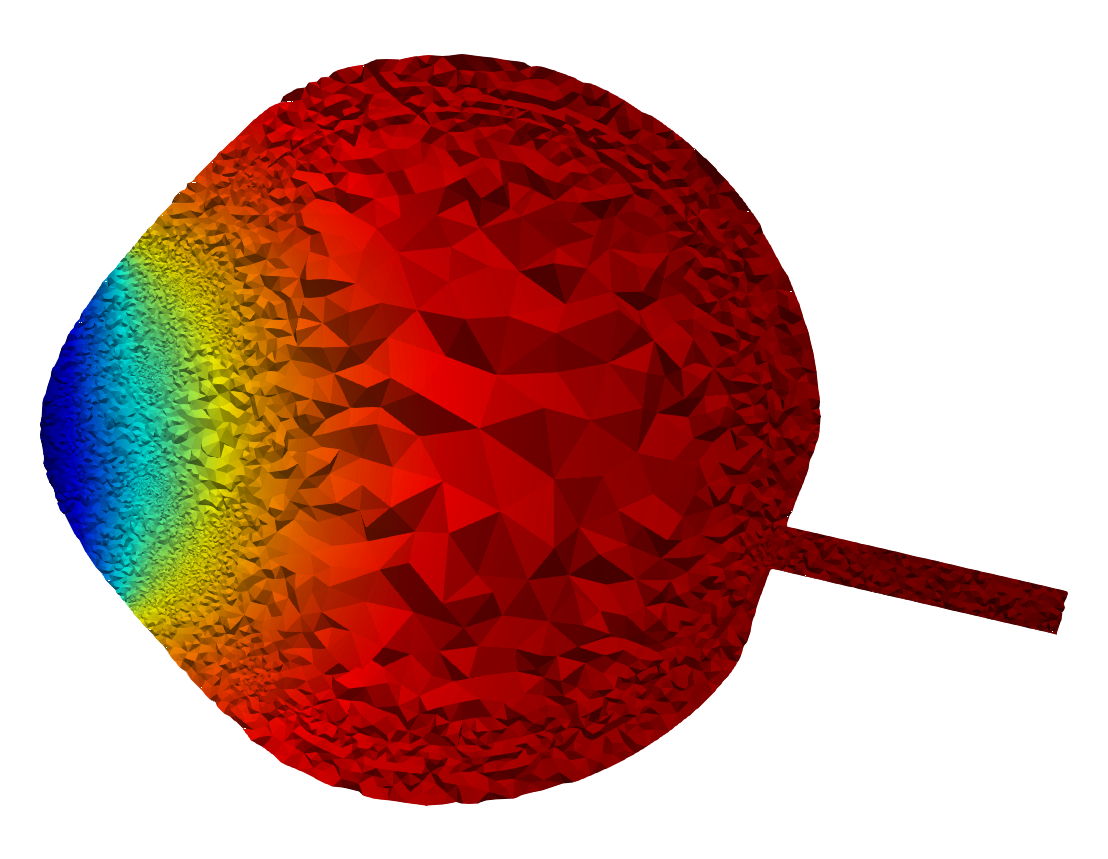};

        \draw[->] (0.8, 0.9) -- (0.8, 0.8) node[midway, anchor=west] {$\vct{g}$};

        \draw[->, red] (0.1, 0.1) -- (0.2, 0.1) node[pos=1, anchor=north] {$x$};
        \draw[->, green!60!black] (0.1, 0.1) -- (0.1, 0.2) node[pos=1, anchor=east] {$y$};
        \draw[blue] (0.1, 0.1) node {$\odot$};
        \draw[blue] (0.1, 0.1) node[anchor=north east] {$z$};
    \end{axis}
    \end{tikzpicture}
    \caption{Distibution of the computed temperature over the eyeball in the standing position, on a vertical cut. Mesh discretization is also presented.}
    \label{fig:temp}
    \vspace*{-1\baselineskip}
\end{wrapfigure}
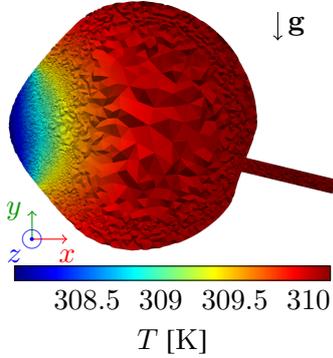

% \vspace{-0.5\baselineskip}
In the spirit of \cite{Abd,Wang}, we make the following assumptions:
the fluid is incompressible, its density is constant, and we neglect the effect of the cornea and the lens on the flow.
% The Boussinesq approximation is used to account for the temperature dependence of the fluid density.
We hence consider that the steady flow of the aqueous humor is governed by the Navier-Stokes equations, coupled with the heat transfer:
\begin{subequations}
    \begin{alignat}{2}
        \rho (\vct{u}\cdot\nabla \vct{u}) - \mu \nabla^2 \vct{u} + \nabla p &= -\rho\beta(T-T_\text{ref})\vct{g} &\quad & \text{in $\Omega_\text{AH}$}, \label{eq:pbStokes:stokes}\\
        \nabla\cdot \vct{u}                         &= 0            &      & \text{in $\Omega_\text{AH}$},\\
        % -2\mu \nabla\cdot(\textstyle\frac12\left(\nabla \vct{u} + \nabla \vct{u}^T\right)) + \nabla p &= -\rho_0\beta(T-T_\text{ref})\vct{g} &\quad & \text{in $\Omega_\text{AH}$}, \label{eq:pbStokes:stokes}\\
        \rho C_p \vct{u}\cdot\nabla T - k\nabla^2 T &= 0            &      & \text{in $\Omega$}, \label{eq:pbStokes:heat}
    \end{alignat}
\label{eq:pbStokes}
\end{subequations}

\vspace*{-1.5\baselineskip}
where $\mu$ [\si{\newton.\second\per\meter^2}] is the dynamic viscosity of the fluid, $\rho$ [\si{\kilogram\per\meter^3}] its density (both at reference temperature $T_\text{ref}$ [\si{\kelvin}]), $C_p$ [\si{\joule\per\kilogram\per\kelvin}] its specific heat, $k$ [\si{\watt\per\meter\per\kelvin}] its thermal conductivity.
The quantity $T$ [\si{\kelvin}] is the temperature of the eye, while $p$ [\si{\pascal}] is the pressure of the aqueous humor fluid (also expressed in \si{\mmHg} in biologic contexts),
and $\vct{u}$ [\si{\meter\per\second}] is its velocity.
The right-hand side term in Eq. \eqref{eq:pbStokes:stokes} represents the gravitational force per unit volume, along with the Boussinesq approximation,
utilized to account for the buoyancy effects.
The coefficient $\beta$ [\si{\kelvin^{-1}}] is the fluid volume expansion coefficient, and $\vct{g}$ [\si{\meter.\second^{-2}}] the gravitational acceleration vector.
Depending on the position of the patient (standing, laying supine or prone respectively), $\vct{g}$ can be either vertical ($\vct{g}=\left[0, -g, 0\right]^T$) or horizontal ($\vct{g}=\left[g, 0, 0\right]^T$, $\vct{g}=\left[-g, 0, 0\right]^T$ respectively), where $g$ is the gravitational acceleration.

We impose the following boundary conditions for the fluid velocity:
% \begin{subequations}
% \begin{alignat}{2}
%     \vct{u} &= \vct{0}      &\quad & \text{on }\Gamma_C\cup \Gamma_I\cup \Gamma_L, \label{eq:BcUp:BCu}\\
%     \vct{u} \times \vct{n}  &= \vct{0} & & \text{on }\Gamma_\text{VH}\cup\Gamma_\text{Sc}, \label{eq:BcUp:BCu2}\\
%     p       &= p_\text{in}  &      & \text{on }\Gamma_\text{VH},\quad p = p_\text{out} \text{ on }\Gamma_\text{Sc} \label{eq:BcUp:BCp}.
% \end{alignat}
% \begin{alignat}{4}
%     \vct{u} &= \vct{0}      & \quad & \text{on }\Gamma_C\cup \Gamma_I\cup \Gamma_L,\quad & \vct{u} \times \vct{n}  &= \vct{0}      & \quad & \text{on }\Gamma_\text{VH}\cup\Gamma_\text{Sc},\label{eq:BcUp:BCu}\\
%     p       &= p_\text{in}  &       & \text{on }\Gamma_\text{VH},                        & p                       &= p_\text{out} & \quad & \text{on }\Gamma_\text{Sc}. \label{eq:BcUp:BCp}
%     \textbf{mettre ce qu'on a fait vraiment}
% \end{alignat}
\begin{alignat}{2}
    \vct{u} &= \vct{0}      &\quad & \text{on }\Gamma_C\cup \Gamma_I\cup \Gamma_L\cup\Gamma_\text{VH}\cup\Gamma_\text{Sc},
    \label{eq:BcUp}
\end{alignat}
% \end{subequations}

\vspace*{-.5\baselineskip}
% The flow of aqueous humor is driven on the one hand by the pressure difference between the anterior and posterior chambers of the eye, provided by the ciliary body, as shown in Eq. \eqref{eq:BcUp:BCp},
% and on the other hand by the temperature difference.
% The inlet pressure imposed on $\Gamma_\text{VH}$, $p_\text{in}$ [\si{\kilogram.\meter^{-1}.\second^{-2}}] (also expressed in \si{\mmHg} in biologic contexts), is the pressure of the aqueous humor in the anterior chamber of the eye, namely the IOP.
that, together with Eq.~\eqref{eq:pbStokes:stokes}, only incorporates thermal and postural mechanisms driving the AH flow.
The hydraulic pressure difference due to production and drainage of AH is not currently accounted for in the model, since previous works \cite{Ooi, kumar} already pointed out that buoyancy is observed to be the dominant source of the convective motion, regardless of the postural orientation.

Moreover, the heat transfer inside the eye is governed by the heat equation \eqref{eq:pbStokes:heat}, with the following boundary conditions,
taking into account the convective heat transfer with the surrounding tissues, over a domain denoted as $\Gamma_\text{body}$; and the heat production due to the metabolism of the eye,
as well as the heat transfer with the ambient air, over a domain called $\Gamma_\text{amb}$:
\begin{subequations}
\begin{alignat}{2}
    -k\textstyle\frac{\partial T}{\partial \vct{n}} &= h_\text{bl}(T-T_\text{bl})      &\quad & \text{on }\Gamma_\text{body}, \label{eq:BcT:BCTbody}\\
    -k\textstyle\frac{\partial T}{\partial \vct{n}} &= (h_\text{amb}+h_\text{r})(T-T_\text{bl}) + E     &\quad & \text{on }\Gamma_\text{amb}. \label{eq:BcT:BCTamb}
\end{alignat}
\label{eq:BcT}
\end{subequations}

\vspace*{-1.5\baselineskip}
Thanks to the open source library Feel++~\cite{Fpp}, we develop a computational framework allowing to simulate the model (\ref{eq:pbStokes}--\ref{eq:BcUp}--\ref{eq:BcT}), using the finite element method.
We solve the monolithic systems with a Newton method, and
% At each nonlinear-iteration we solve for the system on the Jacobian and use the Feel++ based PETSc to solve and precondition the system.
% We use the fieldsplit strategy from PETSc in parallel and develop a multi-physic preconditioner based on:
use a preconditioner based on the PDE system proposed by~\cite[Ch.~11]{elman}.
The overall preconditioner is block diagonal with respect to the fluid and heat unknowns.
The fluid part preconditioner is block upper triangular with the pressure block being the Schur complement~\cite{elman}, the other ones corresponding to the full velocity preconditioned by an additive Schwarz method (ASM) and the velocity pressure coupling respectively
 The heat block is preconditioned also with an ASM.

The geometry presented in Fig.~\ref{fig:geo-eye} is discretized, with a significant effort made on the AC and PC domains, where the coupled fluid-thermal problem is solved.
This results in a mesh composed of \pgfmathprintnumber{5781727} elements.
The refinement of the mesh is presented among results in Fig.~\ref{fig:temp}.

% The linear algebraic systems at each non-linear iteration reads
% \begin{equation}
%     \label{eq:system}
%     \mathcal{J} \left\{
%     \begin{bmatrix}
%         A^{n}_\nu & B^T & D^{n} \\
%         B           & 0   & 0                              \\
%         E^{n}           & 0   & F^{n}
%     \end{bmatrix}\right.
%     \begin{bmatrix}
%         \delta_u^{n+1} \\
%         \delta_p^{n+1} \\
%         \delta_T^{n+1}
%     \end{bmatrix}
%     =
%     \begin{bmatrix}
%         f^{n} \\
%         0       \\
%         f^{n}
%     \end{bmatrix}
%     .
% \end{equation}
%
%%The pressure boundary conditions are imposed using the method presented in \cite{Pressure} \emph{peut être à enlever}.
\begin{wraptable}{R}{0pt}
    \centering
    \begin{tabular}{cc}
        \textbf{Parameter} & \textbf{Value} \\
        \hline
        $\mu$ & \SI{0.001}{\kilogram.\meter^{-1}.\second^{-1}} \\
        $\rho$ & \SI{1000}{\kilogram.\meter^{-3}} \\
        $C_p$ & \SI{4178}{\joule\per\kilogram\per\kelvin} \\
        $\beta$ & \SI{3e-4}{\kelvin^{-1}} \\
        $k_\text{AH}$ & \SI{0.576}{\watt\per\meter\per\kelvin} \\
        $g$ & \SI{9.81}{\meter.\second^{-2}} \\
        $T_\text{ref}$ & \SI{298}{\kelvin} \\
        $h_\text{bl}$ & \SI{65}{\watt.\meter^{-2}.\kelvin^{-1}} \\
        $h_\text{amb}$ & \SI{10}{\watt.\meter^{-2}.\kelvin^{-1}} \\
        $h_\text{r}$ & \SI{6}{\watt.\meter^{-2}.\kelvin^{-1}} \\
        $E$ & \SI{40}{\watt.\meter^{-3}} \\
        $T_\text{bl}$ & \SI{310}{\kelvin} \\
        $T_\text{amb}$ & \SI{307}{\kelvin} \\
        % $p_\text{in}$ & \SI{15.5}{\mmHg} \\
        % $p_\text{out}$ & \SI{15}{\mmHg} \\
        % $p_\text{out}$ & $\SI{13.5}{\mmHg} = \SI{1799.85}{\kilogram\per\meter\per\second^2}$  \\
    \end{tabular}
    \caption{Parameters used for the simulations. Sources can be found in \cite{Ooi,Wang}.}
    \label{tab:parameters}
    \vspace{-2\baselineskip}
\end{wraptable}

\section{Results and conclusions}

We present the results of the simulations, for a set of parameters presented in Tab.~\ref{tab:parameters}.
These parameters are chosen to fit the case simulating a healthy patient, with nominal physiological values.

% To vizualize the impact of the aquehous humor flow on the heat transfer inside the eye, we plot the temperature computed over an horizontal line passing through the center of the eye.
% We compute on the one hand the model (\ref{eq:pbStokes:heat}--\ref{eq:BcT}) simulating the sole heat transfer,
% and on the other hand the model (\ref{eq:pbStokes}--\ref{eq:BcUp}--\ref{eq:BcT}) simulating the coupled model, for two different positions of the patient, laying prone and supine.
% Results are presented in Fig.~\ref{fig:temp}.
% It shows that is has an impact on the distribution of the temperature especially at the front part of the eye.
% In \cite{Ooi} where the flow model is considered in the sole AC, the authors have found similar results.

We show in Fig.~\ref{fig:temp} the distribution of the computed temperature over the eyeball in the standing position, on a vertical cut.
The temperature is higher at the back part of the eye, as this part is inside the human body, while the front part, where heat exchanges with the ambient air are present, is colder.
This is consistent with previous findings \cite{Heat,Ooi}.

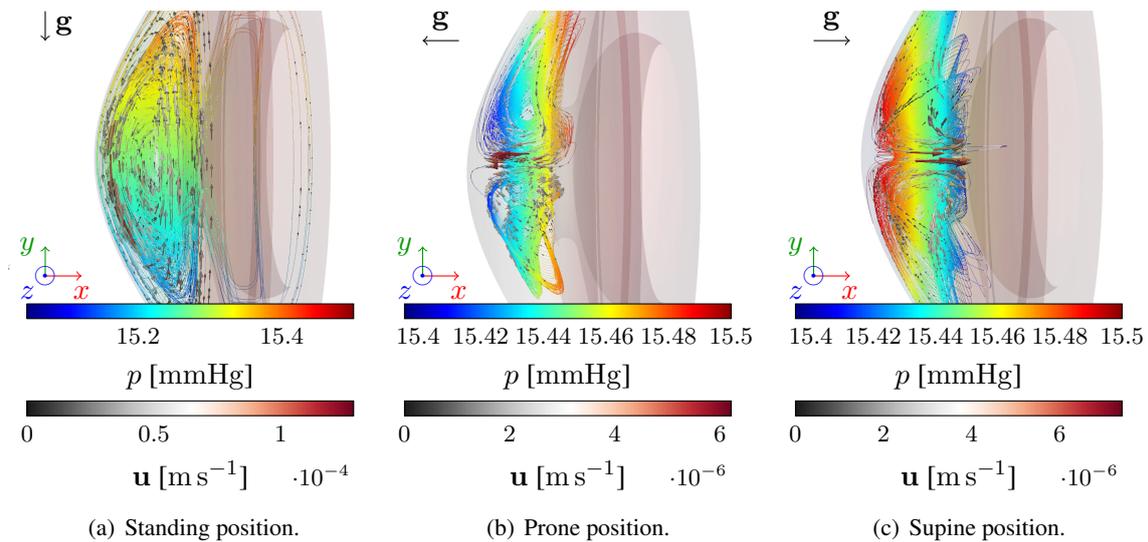
\begin{figure}
    \def\subfigwidth{0.41\textwidth}
    \centering
    \subfigure[Standing position.\label{fig:res:standing}]{%
    \begin{tikzpicture}
    \begin{axis}[
        colorbar,
        colormap/jet, % Choose the colormap you prefer
        % axis equal image,
        enlargelimits=false,
        colorbar horizontal,
        point meta max=15.5,
        point meta min=15.045621357,
        axis line style = {draw=none},
        tick style = {draw=none},
        xtick = \empty, ytick = \empty,
        colorbar style={
            % xlabel style={
            %     at={(0.5,1.1)},
            %     anchor=south,
            % },
            xlabel = {$p$ [\si{\mmHg}]},
            height=0.05*\pgfkeysvalueof{/pgfplots/parent axis height},
            width=0.9*\pgfkeysvalueof{/pgfplots/parent axis width},
            at={(0.5,-0.02)},
            anchor=center,
            tick label style={font=\footnotesize},
        },
        colorbar/draw/.append code={
            \begin{axis}[
                colormap={Gray and Red}{
                    rgb255(-1cm)=(26,26,26);
                    rgb255(-0.87451cm)=(58,58,58);
                    rgb255(-0.74902cm)=(91,91,91);
                    rgb255(-0.623529cm)=(128,128,128);
                    rgb255(-0.498039cm)=(161,161,161);
                    rgb255(-0.372549cm)=(191,191,191);
                    rgb255(-0.247059cm)=(215,215,215);
                    rgb255(-0.121569cm)=(236,236,236);
                    rgb255(0.00392157cm)=(254,254,253);
                    rgb255(0.129412cm)=(253,231,218);
                    rgb255(0.254902cm)=(250,204,180);
                    rgb255(0.380392cm)=(244,170,136);
                    rgb255(0.505882cm)=(228,128,101);
                    rgb255(0.631373cm)=(208,84,71);
                    rgb255(0.756863cm)=(185,39,50);
                    rgb255(0.882353cm)=(147,14,38);
                    rgb255(1cm)=(103,0,31);
                },
                colorbar horizontal,
                point meta min=1.023663860179765e-8,
                point meta max=0.00012864508759841435,
                every colorbar,
                anchor=center,
                colorbar shift,
                colorbar=false,
                xlabel = {$\vct{u}$ [\si{\meter\per\second}]},
                height=0.05*\pgfkeysvalueof{/pgfplots/parent axis height},
                width=0.9*\pgfkeysvalueof{/pgfplots/parent axis width},
                at={(0.5,-0.1*3.5*\pgfkeysvalueof{/pgfplots/parent axis height})},
                tick label style={font=\footnotesize},
            ]
              \pgfkeysvalueof{/pgfplots/colorbar addplot}
            \end{axis}
          },
        width=\subfigwidth
    ]
        \addplot graphics [includegraphics cmd=\pgfimage, xmin=0, xmax=1, ymin=0, ymax=1] {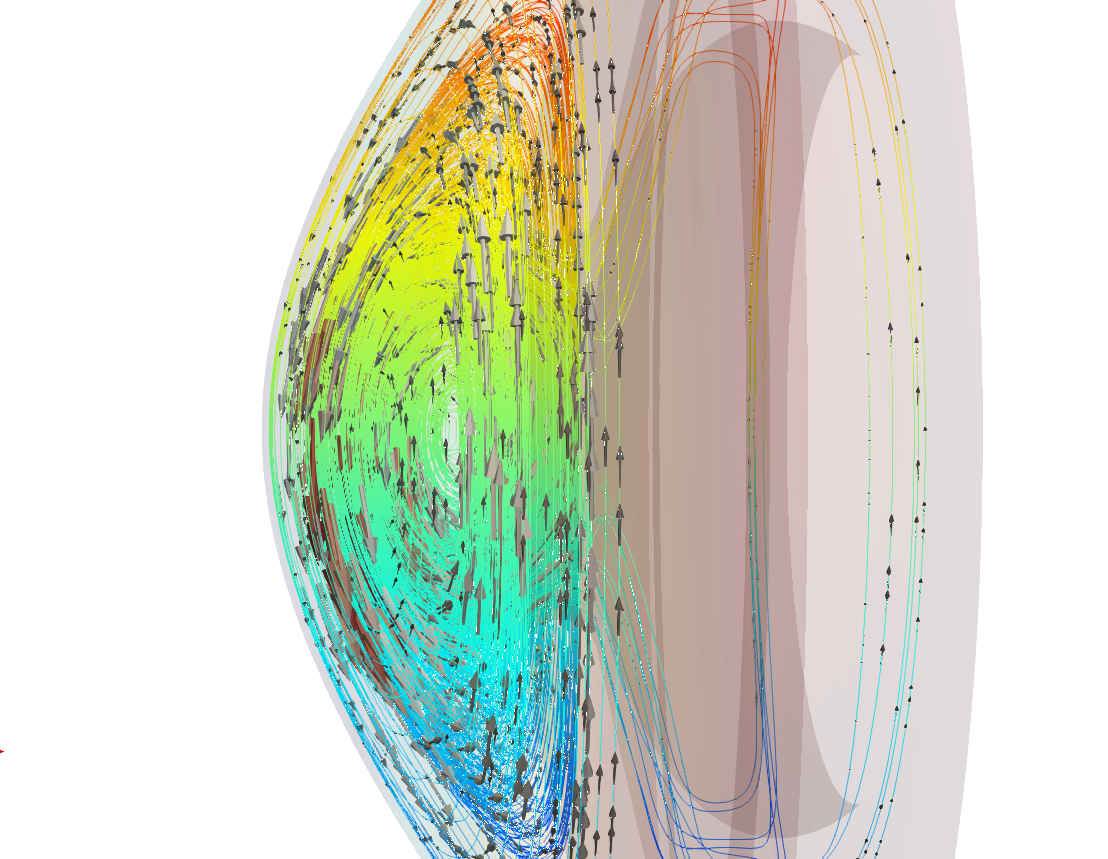};
        \draw[->] (0.1, 1) -- (0.1, 0.9) node[midway, anchor=west] {$\vct{g}$};

        \draw[->, red] (0.1, 0.1) -- (0.2, 0.1) node[pos=1, anchor=north] {$x$};
        \draw[->, green!60!black] (0.1, 0.1) -- (0.1, 0.2) node[pos=1, anchor=east] {$y$};
        \draw[blue] (0.1, 0.1) node {$\odot$};
        \draw[blue] (0.1, 0.1) node[anchor=north east] {$z$};
    \end{axis}
    \end{tikzpicture}
    }
    \subfigure[Prone position.\label{fig:res:prone}]{%
    \begin{tikzpicture}
    \begin{axis}[
        colorbar,
        colormap/jet, % Choose the colormap you prefer
        % axis equal image,
        enlargelimits=false,
        colorbar horizontal,
        point meta max=15.5,
        point meta min=15.395625664,
        axis line style = {draw=none},
        tick style = {draw=none},
        xtick = \empty, ytick = \empty,
        colorbar style={
            % xlabel style={
            %     at={(0.5,1.1)},
            %     anchor=south,
            % },
            xlabel = {$p$ [\si{\mmHg}]},
            height=0.05*\pgfkeysvalueof{/pgfplots/parent axis height},
            width=0.9*\pgfkeysvalueof{/pgfplots/parent axis width},
            at={(0.5,-0.02)},
            anchor=center,
            tick label style={font=\footnotesize},
        },
        colorbar/draw/.append code={
            \begin{axis}[
                colormap={Gray and Red}{
                    rgb255(-1cm)=(26,26,26);
                    rgb255(-0.87451cm)=(58,58,58);
                    rgb255(-0.74902cm)=(91,91,91);
                    rgb255(-0.623529cm)=(128,128,128);
                    rgb255(-0.498039cm)=(161,161,161);
                    rgb255(-0.372549cm)=(191,191,191);
                    rgb255(-0.247059cm)=(215,215,215);
                    rgb255(-0.121569cm)=(236,236,236);
                    rgb255(0.00392157cm)=(254,254,253);
                    rgb255(0.129412cm)=(253,231,218);
                    rgb255(0.254902cm)=(250,204,180);
                    rgb255(0.380392cm)=(244,170,136);
                    rgb255(0.505882cm)=(228,128,101);
                    rgb255(0.631373cm)=(208,84,71);
                    rgb255(0.756863cm)=(185,39,50);
                    rgb255(0.882353cm)=(147,14,38);
                    rgb255(1cm)=(103,0,31);
                },
                colorbar horizontal,
                point meta min=0,
                point meta max=0.0000062121127024231835,
                every colorbar,
                anchor=center,
                colorbar shift,
                colorbar=false,
                xlabel = {$\vct{u}$ [\si{\meter\per\second}]},
                height=0.05*\pgfkeysvalueof{/pgfplots/parent axis height},
                width=0.9*\pgfkeysvalueof{/pgfplots/parent axis width},
                at={(0.5,-0.1*3.5*\pgfkeysvalueof{/pgfplots/parent axis height})},
                tick label style={font=\footnotesize},
            ]
                \pgfkeysvalueof{/pgfplots/colorbar addplot}
            \end{axis}
            },
        width=\subfigwidth
    ]
        \addplot graphics [includegraphics cmd=\pgfimage, xmin=0, xmax=1, ymin=0, ymax=1] {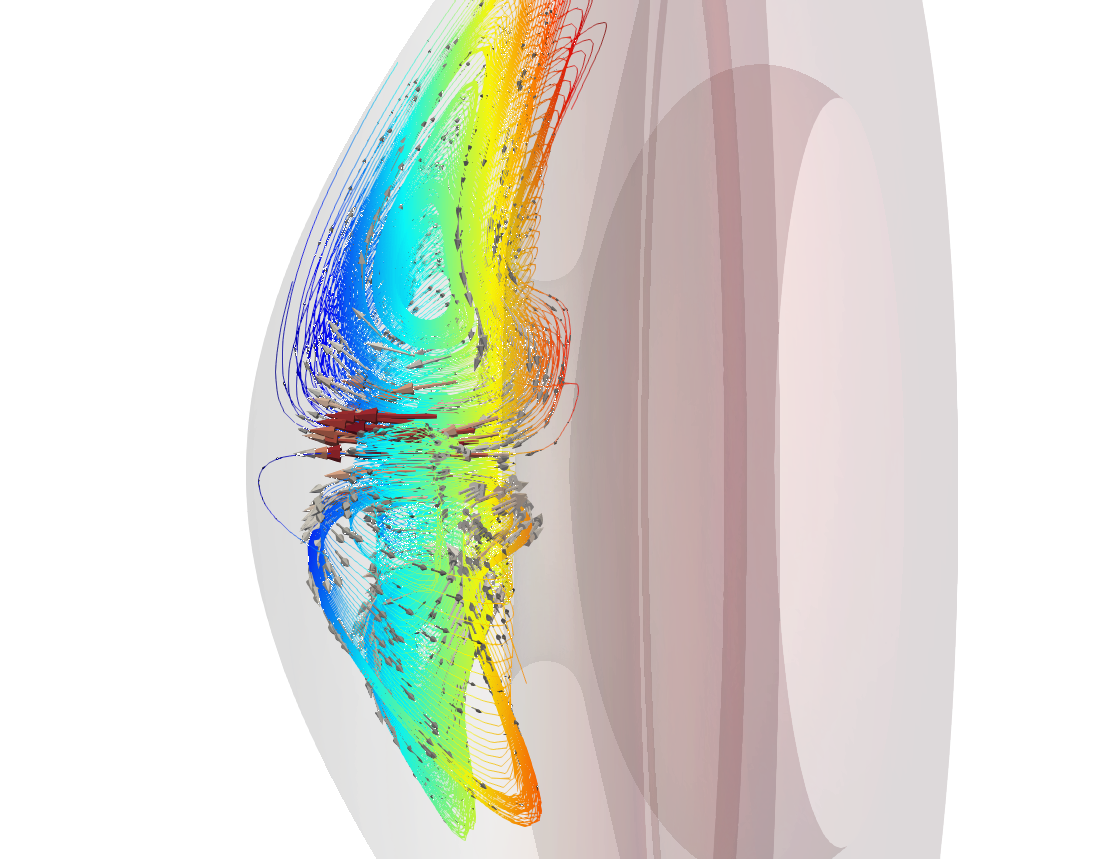};
        \draw[->] (0.2, 0.9) -- (0.1, 0.9) node[midway, anchor=south] {$\vct{g}$};

        \draw[->, red] (0.1, 0.1) -- (0.2, 0.1) node[pos=1, anchor=north] {$x$};
        \draw[->, green!60!black] (0.1, 0.1) -- (0.1, 0.2) node[pos=1, anchor=east] {$y$};
        \draw[blue] (0.1, 0.1) node {$\odot$};
        \draw[blue] (0.1, 0.1) node[anchor=north east] {$z$};
    \end{axis}
    \end{tikzpicture}
    }
    \subfigure[Supine position.\label{fig:res:supine}]{%
    \begin{tikzpicture}
    \begin{axis}[
        colorbar,
        colormap/jet, % Choose the colormap you prefer
        % axis equal image,
        enlargelimits=false,
        colorbar horizontal,
        point meta max=15.5,
        point meta min=15.395051958,
        axis line style = {draw=none},
        tick style = {draw=none},
        xtick = \empty, ytick = \empty,
        colorbar style={
            % xlabel style={
            %     at={(0.5,1.1)},
            %     anchor=south,
            % },
            xlabel = {$p$ [\si{\mmHg}]},
            height=0.05*\pgfkeysvalueof{/pgfplots/parent axis height},
            width=0.9*\pgfkeysvalueof{/pgfplots/parent axis width},
            at={(0.5,-0.02)},
            anchor=center,
            tick label style={font=\footnotesize},
        },
        colorbar/draw/.append code={
            \begin{axis}[
                colormap={Gray and Red}{
                    rgb255(-1cm)=(26,26,26);
                    rgb255(-0.87451cm)=(58,58,58);
                    rgb255(-0.74902cm)=(91,91,91);
                    rgb255(-0.623529cm)=(128,128,128);
                    rgb255(-0.498039cm)=(161,161,161);
                    rgb255(-0.372549cm)=(191,191,191);
                    rgb255(-0.247059cm)=(215,215,215);
                    rgb255(-0.121569cm)=(236,236,236);
                    rgb255(0.00392157cm)=(254,254,253);
                    rgb255(0.129412cm)=(253,231,218);
                    rgb255(0.254902cm)=(250,204,180);
                    rgb255(0.380392cm)=(244,170,136);
                    rgb255(0.505882cm)=(228,128,101);
                    rgb255(0.631373cm)=(208,84,71);
                    rgb255(0.756863cm)=(185,39,50);
                    rgb255(0.882353cm)=(147,14,38);
                    rgb255(1cm)=(103,0,31);
                },
                colorbar horizontal,
                point meta min=1.6125707421871079e-9,
                point meta max=0.000007417459292689545,
                every colorbar,
                anchor=center,
                colorbar shift,
                colorbar=false,
                xlabel = {$\vct{u}$ [\si{\meter\per\second}]},
                height=0.05*\pgfkeysvalueof{/pgfplots/parent axis height},
                width=0.9*\pgfkeysvalueof{/pgfplots/parent axis width},
                at={(0.5,-0.1*3.5*\pgfkeysvalueof{/pgfplots/parent axis height})},
                tick label style={font=\footnotesize},
            ]
                \pgfkeysvalueof{/pgfplots/colorbar addplot}
            \end{axis}
            },
        width=\subfigwidth
    ]
        \addplot graphics [includegraphics cmd=\pgfimage, xmin=0, xmax=1, ymin=0, ymax=1] {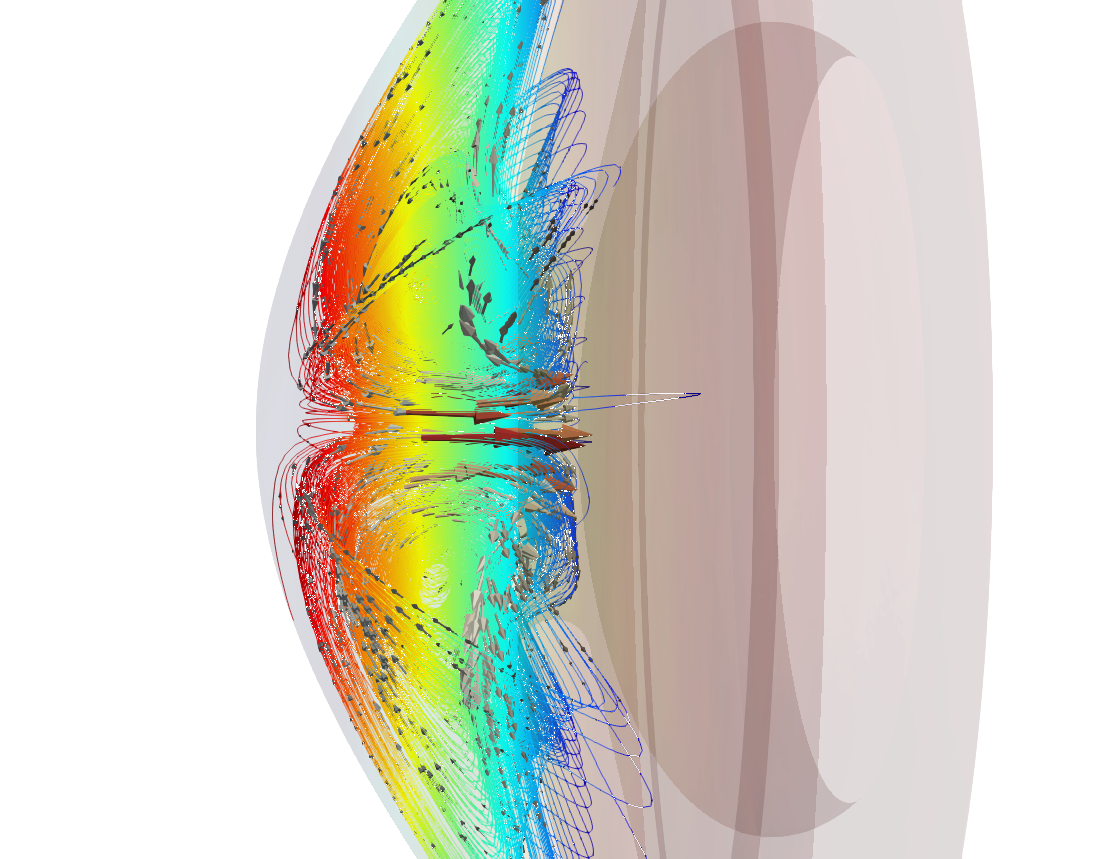};
        \draw[->] (0.1, 0.9) -- (0.2, 0.9) node[midway, anchor=south] {$\vct{g}$};

        \draw[->, red] (0.1, 0.1) -- (0.2, 0.1) node[pos=1, anchor=north] {$x$};
        \draw[->, green!60!black] (0.1, 0.1) -- (0.1, 0.2) node[pos=1, anchor=east] {$y$};
        \draw[blue] (0.1, 0.1) node {$\odot$};
        \draw[blue] (0.1, 0.1) node[anchor=north east] {$z$};
    \end{axis}
    \end{tikzpicture}
    }
    \caption{Results of simulation for various orientation of the eye. Streamlines are colored according to the pressure, and the arrows show the fluid velocity.}
    \label{fig:res}
\end{figure}

Now we focus on the impact of the position of the subject on the flow of the AH in the AC.
We present in Fig.~\ref{fig:res}, the results for standing and laying positions, prone and supine.
The main striking result is the difference in the flow patterns: it follows the gravitational force, and is more important in the standing position (Fig.~\ref{fig:res:standing}).
Moreover, we notice the formation of Krukenberg’s spindle and phenomena of recirculation in the AC, as observed in the literature \cite{Abd,Wang,Murgoitio}.

\paragraph*{Conclusions and perspectives.}
This work presented a modeling and computational framework allowing to simulate the heat transfer inside the human eyeball, coupled with  the flow of the AH, in the anterior and posterior chambers of a healthy eye.
The results on flow patterns and temperature distribution are in line with previous numerical studies.
In addition, the impact of the postural orientation on the flow recirculations is accurately captured.
One main drawback of the present model is the computational cost, which is too high for real-time simulations.
In order to overcome this issue, we are currently working on the development of a model order reduction technique for the present problem, to enable real-time simulations of coupled flow-heat transfer inside the human eyeball.
From the clinical perspective, such an approach could be utilized to assess the effect of the topical administration of different drugs (eye drops) that are nowadays the standards therapeutical approach for ocular diseases such as glaucoma.

\paragraph*{Acknowledgements}
The authors would like to acknowledge the support of the platform Cemosis at University of Strasbourg and the French Ministry of Higher Education, Research and Innovation.

% REFERENCES: managed using thebibliography with cite_key and \bibitem{cite_key} command.
% references are arranged by order of appearance in the text.
% Please follow the format of each entries.

\end{document}

%% file: fig/eye.pdf_tex
%% Creator: Inkscape 1.3.1 (9b9bdc1480, 2023-11-25, custom), www.inkscape.org
%% PDF/EPS/PS + LaTeX output extension by Johan Engelen, 2010
%% Accompanies image file 'eye.pdf' (pdf, eps, ps)
%%
%% To include the image in your LaTeX document, write
%%   \input{<filename>.pdf_tex}
%%  instead of
%%   \includegraphics{<filename>.pdf}
%% To scale the image, write
%%   \def\svgwidth{<desired width>}
%%   \input{<filename>.pdf_tex}
%%  instead of
%%   \includegraphics[width=<desired width>]{<filename>.pdf}
%%
%% Images with a different path to the parent latex file can
%% be accessed with the `import' package (which may need to be
%% installed) using
%%   \usepackage{import}
%% in the preamble, and then including the image with
%%   \import{<path to file>}{<filename>.pdf_tex}
%% Alternatively, one can specify
%%   \graphicspath{{<path to file>/}}
%%
%% For more information, please see info/svg-inkscape on CTAN:
%%   http://tug.ctan.org/tex-archive/info/svg-inkscape
%%
\begingroup%
  \makeatletter%
  \providecommand\color[2][]{%
    \errmessage{(Inkscape) Color is used for the text in Inkscape, but the package 'color.sty' is not loaded}%
    \renewcommand\color[2][]{}%
  }%
  \providecommand\transparent[1]{%
    \errmessage{(Inkscape) Transparency is used (non-zero) for the text in Inkscape, but the package 'transparent.sty' is not loaded}%
    \renewcommand\transparent[1]{}%
  }%
  \providecommand\rotatebox[2]{#2}%
  \newcommand*\fsize{\dimexpr\f@size pt\relax}%
  \newcommand*\lineheight[1]{\fontsize{\fsize}{#1\fsize}\selectfont}%
  \ifx\svgwidth\undefined%
    \setlength{\unitlength}{778.63751221bp}%
    \ifx\svgscale\undefined%
      \relax%
    \else%
      \setlength{\unitlength}{\unitlength * \real{\svgscale}}%
    \fi%
  \else%
    \setlength{\unitlength}{\svgwidth}%
  \fi%
  \global\let\svgwidth\undefined%
  \global\let\svgscale\undefined%
  \makeatother%
  \begin{picture}(1,0.59082573)%
    \lineheight{1}%
    \setlength\tabcolsep{0pt}%
    \put(0,0){\includegraphics[width=\unitlength,page=1]{eye.pdf}}%
    \put(0,0){\includegraphics[width=\unitlength,page=2]{eye.pdf}}%
    \put(0.03281095,0.30790147){\color[rgb]{0,0,0}\makebox(0,0)[t]{\lineheight{1.25}\smash{\begin{tabular}[t]{c}Ciliary\\body\end{tabular}}}}%
    \put(0.68975632,0.38259866){\color[rgb]{0,0,0}\makebox(0,0)[lt]{\lineheight{1.25}\smash{\begin{tabular}[t]{l}Lamina Cribrosa\end{tabular}}}}%
    \put(0,0){\includegraphics[width=\unitlength,page=3]{eye.pdf}}%
    \put(0.75191582,0.3130218){\color[rgb]{0,0,0}\makebox(0,0)[lt]{\lineheight{1.25}\smash{\begin{tabular}[t]{l}Optic Nerve\end{tabular}}}}%
    \put(0.43676274,0.36714171){\color[rgb]{0,0,0}\makebox(0,0)[t]{\lineheight{1.25}\smash{\begin{tabular}[t]{c}Vitreous\\humor\end{tabular}}}}%
    \put(0,0){\includegraphics[width=\unitlength,page=4]{eye.pdf}}%
    \put(0.34401795,0.19743467){\color[rgb]{0,0,0}\makebox(0,0)[lt]{\lineheight{1.25}\smash{\begin{tabular}[t]{l}Lens\end{tabular}}}}%
    \put(0,0){\includegraphics[width=\unitlength,page=5]{eye.pdf}}%
    \put(0.66026044,0.53175482){\color[rgb]{0,0,0}\makebox(0,0)[lt]{\lineheight{1.25}\smash{\begin{tabular}[t]{l}Choroid\end{tabular}}}}%
    \put(0,0){\includegraphics[width=\unitlength,page=6]{eye.pdf}}%
    \put(0.69906994,0.05141589){\color[rgb]{0,0,0}\makebox(0,0)[lt]{\lineheight{1.25}\smash{\begin{tabular}[t]{l}Retina\end{tabular}}}}%
    \put(0.0556473,0.53499595){\color[rgb]{0,0,0}\makebox(0,0)[lt]{\lineheight{1.25}\smash{\begin{tabular}[t]{l}Sclera\end{tabular}}}}%
    \put(0,0){\includegraphics[width=\unitlength,page=7]{eye.pdf}}%
    \put(0.13181095,0.06790147){\color[rgb]{0,0,0}\makebox(0,0)[t]{\lineheight{1.25}\smash{\begin{tabular}[t]{c}Aqueous\\humor\end{tabular}}}}%
    \put(0,0){\includegraphics[width=\unitlength,page=8]{eye.pdf}}%
    \put(-0.00216348,0.43464768){\color[rgb]{0,0,0}\makebox(0,0)[lt]{\lineheight{1.25}\smash{\begin{tabular}[t]{l}Cornea\end{tabular}}}}%
    \put(0,0){\includegraphics[width=\unitlength,page=9]{eye.pdf}}%
    \put(0.94684963,0.04818223){\color[rgb]{1,0,0}\makebox(0,0)[lt]{\lineheight{1.25}\smash{\begin{tabular}[t]{l}$x$\end{tabular}}}}%
    \put(0.85384689,0.13300943){\color[rgb]{0,0.7254902,0}\makebox(0,0)[lt]{\lineheight{1.25}\smash{\begin{tabular}[t]{l}$y$\end{tabular}}}}%
    \put(0.88882513,-0.02){\color[rgb]{0,0,1}\makebox(0,0)[lt]{\lineheight{1.25}\smash{\begin{tabular}[t]{l}$z$\end{tabular}}}}%
    \put(0,0){\includegraphics[width=\unitlength,page=10]{eye.pdf}}%
  \end{picture}%
\endgroup%

%% file: fig/aqueous_humor.pdf_tex
%% Creator: Inkscape 1.3.1 (9b9bdc1480, 2023-11-25, custom), www.inkscape.org
%% PDF/EPS/PS + LaTeX output extension by Johan Engelen, 2010
%% Accompanies image file 'aqueous_humor.pdf' (pdf, eps, ps)
%%
%% To include the image in your LaTeX document, write
%%   \input{<filename>.pdf_tex}
%%  instead of
%%   \includegraphics{<filename>.pdf}
%% To scale the image, write
%%   \def\svgwidth{<desired width>}
%%   \input{<filename>.pdf_tex}
%%  instead of
%%   \includegraphics[width=<desired width>]{<filename>.pdf}
%%
%% Images with a different path to the parent latex file can
%% be accessed with the `import' package (which may need to be
%% installed) using
%%   \usepackage{import}
%% in the preamble, and then including the image with
%%   \import{<path to file>}{<filename>.pdf_tex}
%% Alternatively, one can specify
%%   \graphicspath{{<path to file>/}}
%%
%% For more information, please see info/svg-inkscape on CTAN:
%%   http://tug.ctan.org/tex-archive/info/svg-inkscape
%%
\begingroup%
  \makeatletter%
  \providecommand\color[2][]{%
    \errmessage{(Inkscape) Color is used for the text in Inkscape, but the package 'color.sty' is not loaded}%
    \renewcommand\color[2][]{}%
  }%
  \providecommand\transparent[1]{%
    \errmessage{(Inkscape) Transparency is used (non-zero) for the text in Inkscape, but the package 'transparent.sty' is not loaded}%
    \renewcommand\transparent[1]{}%
  }%
  \providecommand\rotatebox[2]{#2}%
  \newcommand*\fsize{\dimexpr\f@size pt\relax}%
  \newcommand*\lineheight[1]{\fontsize{\fsize}{#1\fsize}\selectfont}%
  \ifx\svgwidth\undefined%
    \setlength{\unitlength}{386.11268532bp}%
    \ifx\svgscale\undefined%
      \relax%
    \else%
      \setlength{\unitlength}{\unitlength * \real{\svgscale}}%
    \fi%
  \else%
    \setlength{\unitlength}{\svgwidth}%
  \fi%
  \global\let\svgwidth\undefined%
  \global\let\svgscale\undefined%
  \makeatother%
  \begin{picture}(1,1.3461637)%
    \lineheight{1}%
    \setlength\tabcolsep{0pt}%
    \put(0,0){\includegraphics[width=\unitlength,page=1]{aqueous_humor.pdf}}%
    \put(0.03250118,1.15){\makebox(0,0)[lt]{\lineheight{1.25}\smash{\begin{tabular}[t]{l}Anterior\\chamber\end{tabular}}}}%
    \put(0.793955,0.85376655){\makebox(0,0)[lt]{\lineheight{1.25}\smash{\begin{tabular}[t]{l}Posterior\\chamber\end{tabular}}}}%
    \put(0.07,0.68){\color[rgb]{0,0.67058824,0.01568627}\makebox(0,0)[lt]{\lineheight{1.25}\smash{\begin{tabular}[t]{l}$\Gamma_C$\end{tabular}}}}%
    \put(0,0){\includegraphics[width=\unitlength,page=4]{aqueous_humor.pdf}}%
    \put(0.85703755,0.53662101){\color[rgb]{0,0.67058824,0.01568627}\makebox(0,0)[lt]{\lineheight{1.25}\smash{\begin{tabular}[t]{l}$\Gamma_L$\end{tabular}}}}%
    \put(0,0){\includegraphics[width=\unitlength,page=5]{aqueous_humor.pdf}}%
    \put(0.88,1.19){\color[rgb]{0,0.67058824,0.01568627}\makebox(0,0)[lt]{\lineheight{1.25}\smash{\begin{tabular}[t]{l}$\Gamma_I$\\\end{tabular}}}}%
    \put(0.77488034,0.29134354){\color[rgb]{0.86666667,0,0.01568627}\makebox(0,0)[lt]{\lineheight{1.25}\smash{\begin{tabular}[t]{l}$\Gamma_\text{VH}$\end{tabular}}}}%
    \put(0.38,1.28){\color[rgb]{0.01568627,0,0.86666667}\makebox(0,0)[lt]{\lineheight{1.25}\smash{\begin{tabular}[t]{l}$\Gamma_\text{Sc}$\end{tabular}}}}%
    \put(0,0){\includegraphics[width=\unitlength,page=6]{aqueous_humor.pdf}}%
    \put(0.34802232,0.6867851){\makebox(0,0)[lt]{\lineheight{1.25}\smash{\begin{tabular}[t]{l}$\Omega_\text{AH}$\end{tabular}}}}%
    \put(0,0){\includegraphics[width=\unitlength,page=7]{aqueous_humor.pdf}}%
    \put(0.23125884,0.19096018){\color[rgb]{1,0,0}\makebox(0,0)[lt]{\lineheight{1.25}\smash{\begin{tabular}[t]{l}$x$\end{tabular}}}}%
    \put(0.0437089,0.36202338){\color[rgb]{0,0.7254902,0}\makebox(0,0)[lt]{\lineheight{1.25}\smash{\begin{tabular}[t]{l}$y$\end{tabular}}}}%
    \put(0,0){\includegraphics[width=\unitlength,page=8]{aqueous_humor.pdf}}%
    \put(0.03358195,0.0953055){\color[rgb]{0,0,1}\makebox(0,0)[lt]{\lineheight{1.25}\smash{\begin{tabular}[t]{l}$z$\end{tabular}}}}%
    % \put(0,0){\includegraphics[width=\unitlength,page=9]{aqueous_humor.pdf}}%
  \end{picture}%
\endgroup%